%% file: main_final.tex
\documentclass[12pt,draft]{article}
\usepackage{tikzit}
\input{resonator.tikzstyles}
\usepackage[english]{babel}
\usepackage[utf8x]{inputenc} 
\usepackage[T1]{fontenc}
\usepackage{graphicx}
\usepackage{caption, subcaption}
\usepackage[title]{appendix}
\usepackage{pgfplots}
\usepackage{cutwin}
\graphicspath{ {images/} }
\usetikzlibrary{intersections,decorations.pathreplacing,decorations.markings,calc,angles,quotes,arrows.meta,pgfplots.fillbetween,patterns}
\usepackage[a4paper,top=2.9cm,bottom=2cm,left=2.5cm,right=2.5cm,marginparwidth=1.75cm]{geometry}

\usepackage{amsmath}
\usepackage{amsmath,yhmath}
\usepackage{amsthm}
\usepackage{amssymb}
\usepackage{hyperref}
\usepackage{cleveref}
\usepackage{epstopdf}
\usepackage{comment}
\usepackage{todonotes}
\usepackage{color}
\usepackage{float}
\usepackage[sort&compress,numbers]{natbib}
\newcommand{\nm}{\noalign{\smallskip}}
\newcommand{\ds}{\displaystyle}

\newtheorem{thm}{Theorem}

\theoremstyle{definition}
\newtheorem{defn}[thm]{Definition}
\newtheorem{rmk}[thm]{Remark}
\numberwithin{equation}{section}
\numberwithin{thm}{section}

\newcommand{\R}{\mathbb{R}}

\renewcommand{\S}{\mathcal{S}}

\newcommand{\Z}{\mathbb{Z}}
\newcommand{\p}{\partial}
\renewcommand{\epsilon}{\varepsilon}
\newcommand{\dx}{\: \mathrm{d}}

\newcommand{\capacity}{\mathrm{Cap}}

\newcommand{\iu}{\mathrm{i}\mkern1mu}

\makeatletter
\newcommand{\neutralize}[1]{\expandafter\let\csname c@#1\endcsname\count@}
\makeatother

\usepackage{natbib}
\usepackage{graphicx}

\title{On the validity of the tight-binding method for describing systems of subwavelength resonators}

\author{Habib Ammari\thanks{\footnotesize Department of Mathematics, ETH Z\"urich, R\"amistrasse 101, CH-8092 Z\"urich, Switzerland (habib.ammari@math.ethz.ch, ffiorani@student.ethz.ch, erik.orvehed.hiltunen@sam.math.ethz.ch).}\and Francesco Fiorani\footnotemark[1]  \and Erik Orvehed Hiltunen\footnotemark[1]}

\date{}

\begin{document}
	\maketitle
	
\begin{abstract}
	The goal of this paper is to relate the capacitance matrix formalism to the tight-binding approximation. By doing so, we open the way to the use of mathematical techniques and tools from condensed matter theory in the mathematical and numerical analysis of metamaterials, in particular for the understanding of their topological properties. We firstly study how the capacitance matrix formalism, both when the material parameters are static and modulated, can be posed in a Hamiltonian form. Then, we use this result to compare this formalism to the tight-binding approximation. We prove that the correspondence between the capacitance formulation and the tight-binding approximation holds only in the case of dilute resonators. On the other hand, the tight-binding model is often coupled with a nearest-neighbour approximation, whereby long-range interactions are neglected.  Even in the dilute case, we show that long-range interactions between subwavelength resonators are relatively strong and nearest-neighbour approximations are not generally appropriate.
	\end{abstract}

\noindent{\textbf{Mathematics Subject Classification (MSC2000):} 35J05, 35C20, 35P20, 74J20
		
\vspace{0.2cm}
		
	\vspace{0.5cm}
		
\tableofcontents

\section{Introduction}
The study of acoustic wave propagation in metamaterials was first pioneered by Minnaert in 1933. Here, by metamaterials we mean materials with microstructures that act as subwavelength resonators, that is, structures with diameter much smaller than the operating wavelength \cite{lemoult2016soda,yves2017crytalline,phononic1,phononic2}. Due to the subwavelength nature of some of the resonances of their building blocks, these microstructures can produce a high interference for long wavelengths. In \cite{minnaert1933musical}, Minnaert  studied acoustic waves in bubbly media and noticed that the resonant and scattering properties of sound moving inside water is significantly shifted when introducing air bubbles. In the 1950s an artificial screen of bubbles in a lake was used to shield damns from shocks. More recently, sophisticated wave manipulations, such as waveguides, super-resolution,  cloaking and shielding devices were obtained by means of metamaterials \cite{review,review2}.

A fundamental property of subwavelength resonators is that their material parameters differ greatly from the background medium. To investigate their subwavelength resonances, a capacitance matrix formalism describing the properties of finite as well as periodic systems of subwavelength resonators has been introduced \cite{ammari2020biomimetic,ammari2019cochlea,ammari2021functional,ammari2019double}. The capacitance matrix formulation is a {\em discrete approximation} of the wave equation. It has enabled the mathematical justification of many of the remarkable applications of metamaterials listed above \cite{ammari2020exceptional,ammari2020high,ammari2019double,ammari2017effective,ammari2019bloch}, also in the case of time-dependent material parameters \cite{ammari2020time}. See also the recent review paper on the topic \cite{ammari2021functional}.
To achieve this theory, layer-potential techniques and Gohberg-Sigal theory \cite{gohberg1971operator,ammari2018mathematical}  have been used. Floquet-Bloch theory has been also utilised when the material is supposed to possess a periodic  structure. The starting point of the theory is the wave equation. 

Superficially, the capacitance formulation bears resemblance to the tight-binding approximation commonly used in condensed matter theory. Like the capacitance matrix formulation, the tight-binding approximation enables to recast a continuous problem (Schr\"odinger's equation) as a discrete one. It consists in assuming that the main interactions are given by self-interactions. In the case where the electrons are strongly bonded to the atoms, a mode for the whole system can be well approximated by a linear combination of modes relative to single atoms \cite{https://doi.org/10.1002/cpa.21735}. Similarly as for the tight-binding approximation, the capacitance formulation has been used to describe topological properties of subwavelength resonator systems \cite{ammari2020topological,ammari2020robust,ammari2018high,ammari2020honeycomb}. In large structures, the tight-binding approximation is usually coupled with a nearest-neighbour approximation, where long-range interactions are neglected.

The main goal of this paper is to link the capacitance formulation of subwavelength resonators with the tight-binding approximation.  By doing so, we open the way to the use of mathematical techniques and tools from condensed matter theory in the mathematical and numerical analysis of subwavelength metamaterials, in particular for investigating their topological properties. We firstly rewrite the resonant problem for the wave equation in a Hamiltonian form. Based on this, we investigate the validity of the tight-binding approximation. We demonstrate that a tight-binding-type approximation is only valid if the resonators are dilute, \textit{i.e.}, they are significantly smaller than their separating distance. Nevertheless, even in the dilute case, long-range interactions cannot be accurately neglected. 


\section{Setting}
Our goal is to provide a characterisation for the solution to the wave equation in a structure composed of contrasting materials. The material parameters are given by $\rho(x,t)$ and $\kappa(x,t)$. In the example of acoustic waves $\rho$ and $\kappa$ correspond to the density and the bulk modulus. We remark that the equation of study models a range of problems. In particular, it can be used to study both acoustic waves and polarized electromagnetic waves.
	
	We consider the time-dependent wave equation in three dimensions,
	\begin{equation}\label{eq:wave}
		\left(\frac{\p }{\p t } \frac{1}{\kappa(x,t)} \frac{\p}{\p t} - \nabla \cdot \frac{1}{\rho(x,t)} \nabla\right) u(x,t) = 0, \quad x\in \R^3, t\in \R.
	\end{equation}
Here, $\nabla$ denotes the gradient with respect to $x\in \R^3$. We assume that $D\subset \mathbb{R}^d$ is constituted by $N$ disjoint bounded domains $D_i$ for $i=1,...,N$, each $D_i$ being connected and having boundary of Hölder class $\p D_i \in C^{1,s}, 0 < s < 1$. We will assume that the modulation only takes place inside the resonators, which was shown in \cite{ammari2020time} to be an effective way to obtain phenomena not achievable in the static case. The material parameters are hence discontinuous and of the form
\begin{equation} \label{eq:resonatormod_finite}
	\kappa(x,t) = \begin{cases}
		\kappa_0, & x \in \R^3 \setminus \overline{D}, \\  \kappa_r\kappa_i(t), & x\in D_i,
	\end{cases}, \qquad \rho(x,t) = \begin{cases}
		\rho_0, & x \in \R^3 \setminus \overline{D}, \\  \rho_r\rho_i(t), & x\in D_i, \end{cases}
\end{equation}
where $\rho_r$, $\kappa_r$ are positive constants and $\rho_i$, $\kappa_i$, for $1 \leq i \leq N$, are smooth periodic functions with period $\tau$ depending only on time $t$.  The functions $\rho_i(t)$ and $\kappa_i(t)$ correspond to the modulation inside the corresponding resonator $D_i$. We assume that they have finite Fourier series.  It is also worth emphasizing that the fact that the material parameters are assumed to be real is key for the physical meaning of energy gain and loss we will have in deriving a Hamiltonian formulation in the next section.

We seek solutions to \eqref{eq:wave} under the modulation specified in \eqref{eq:resonatormod_finite}. From the regularity of $\rho$ and $\kappa$, we know that $u$ is continuously differentiable in $t$ (see \cite{ammari2020time}). Since $e^{- \iu \omega t}u(x,t) $ is a $\tau$-periodic function of $t$, we have a Fourier series expansion as 
\begin{equation} \label{form} u(x,t)= e^{\iu \omega t}\sum_{n = -\infty}^\infty v_n(x)e^{\iu n\Omega t},\end{equation}
where $\Omega = (2\pi)/\tau$. Throughout this paper, we assume that $\Omega = O(\sqrt{\delta})$. Here, $\omega$ is defined modulo $\Omega$, and is contained in the (time)-Brillouin zone $Y_t^* = \mathbb{C} / \left(\Omega \Z\right)$.

Let $\partial/\partial \nu$ denote the outward normal derivative at $\partial D_i$ and let the subscripts $+$ and $-$ denote taking the limit from outside and inside the boundary $\partial D_i$, respectively. From \eqref{eq:wave} it follows that we have the following transmission conditions on $\p D_i$:
$$\delta \frac{\partial {u}}{\partial \nu} \bigg|_{+} - \frac{1}{\rho_i(t)}\frac{\partial {u}}{\partial \nu} \bigg|_{-} = 0, \qquad x\in \p D_i, \ t\in \R,$$
where $0< \delta := \frac{\rho_r}{\rho_0} \ll 1$ is the asymptotic contrast parameter, which is independent of $t$. The assumption $\delta \ll 1$ is  key to achieve subwavelength quasifrequencies. 

\begin{defn}[subwavelength quasifrequencies]
Let $Y_t^*:= \mathbb{C}/\Omega \mathbb{Z}$. If $\omega = O(\sqrt{\delta}) \in Y_t^*$ satisfies
\begin{itemize}
\item[(i)]
 there is a non-zero solution $u$ to 
(\ref{eq:wave}); 
\item[(ii)]
there exists $M(\delta)$ such that $M(\delta) \Omega \rightarrow 0$ as $\delta \rightarrow 0$ and 
$$
\sum_{n=-\infty}^{\infty} \|v_n\|^2_{L^2 (D)} = \sum_{n=-M}^M  \|v_n\|^2_{L^2 (D)} + o(1),
$$
\end{itemize}
where $v_n$ are defined in (\ref{form}), then $\omega$ is called a {\em subwavelength quasifrequency}.
\end{defn}

Note that condition (ii) means that $u$ only contains components which are either in the subwavelength regime (with frequencies of order of $\sqrt{\delta}$), or which are asymptotically small as $\delta \to 0$, see \cite{ammari2020time}. 

Seeking quasiperiodic solutions in $t$, we apply Fourier transform to \eqref{eq:wave}. In the frequency domain, we then have for any $n\in \Z$ the following equation:
\begin{equation} \label{eq:freq}
	\left\{
	\begin{array} {ll}
		\ds \Delta {v_n}+ \frac{\rho_0(\omega+n\Omega)^2}{\kappa_0} {v_n}  = 0 \ \ \ \  \text{in } \mathbb{R}^3 \setminus \overline{D}, \\[0.3em]
		\ds \Delta v_{i,n}^* +\frac{\rho_r(\omega+n\Omega)^2}{\kappa_r} v_{i,n}^{**}  = 0 \ \ \  \text{in } D_i, \\
		\nm
		\ds  {v_n}|_{+} -{v_n}|_{-}  = 0  \quad \quad\text{on } \partial D, \\
		\nm
		\ds  \delta \frac{\partial {v_n}}{\partial \nu} \bigg|_{+} - \frac{\partial v_{i,n}^* }{\partial \nu} \bigg|_{-} = 0 \quad \text{on } \partial D_i, \\[0.3em]
		v_n(x) \text{ satisfies either the outgoing or the incoming radiation condition}.
	\end{array}
	\right.
\end{equation}
Here, $v_{i,n}^*(x)$ and $v_{i,n}^{**}(x)$ are defined through the convolutions
$$v_{i,n}^*(x) = \sum_{m = -\infty}^\infty r_{i,m} v_{n-m}(x), \quad  v_{i,n}^{**}(x) = \frac{1}{\omega+n\Omega}\sum_{m = -\infty}^\infty k_{i,m}\big(\omega+(n-m)\Omega\big)v_{n-m}(x),$$
where $r_{i,m}$ and $k_{i,m}$ are the Fourier series coefficients of $1/\rho_i$ and $1/\kappa_i$, respectively:
$$\frac{1}{\rho_i(t)} = \sum_{n = -M}^M r_{i,n} e^{\iu n \Omega t}, \quad \frac{1}{\kappa_i(t)} = \sum_{n = -M}^{M} k_{i,n} e^{\iu n \Omega t},$$ with $M =M(\delta) \in \mathbb{N}$ satisfying $M= O(\delta ^{-\gamma/2})$ for some $0<\gamma<1$. Again, we remark that one of the key assumptions is to suppose these last two Fourier series to be finite. Moreover, the radiation conditions are given by
\begin{equation} \label{eq:SRC}
	\lim_{|x|\to\infty} |x|\left(\frac{\partial}{\partial |x|}\pm\iu k_n\right)v_n=0, \quad \text{uniformly in all directions } x/|x|,
\end{equation}
where $k_n = \sqrt{\frac{\rho_0}{\kappa_0}}(\omega+n\Omega)$ and where ``$+$'' corresponds to the incoming while ``$-$'' corresponds to the outgoing condition.

Observe that \eqref{eq:freq} consists of coupled Helmholtz equations at frequencies differing by integer multiples of $\Omega$. The coupling of the Helmholtz equations is specified through $\rho_i$ and $\kappa_i$, $1 \leq i \leq N$.

As shown in \cite{ammari2020time}, the functions $v_{i,n}^*(x)$ are asymptotically constant inside the resonators as $\delta \to 0$. In other words, we have $v_{i,n}^*(x) = c_{i,n} + o(\delta)$. This enables us to employ the capacitance matrix formulation to obtain a system of ordinary differential equations (ODEs) satisfied by the coefficients $c_{i,n}$ as the contrast parameter $\delta$ goes to zero. We define 
$$c_i(t) = \sum_{n \in \Z} c_{i,n}e^{\iu n \Omega t}.$$

In order to introduce the capacitance formulation, we will use the single layer potential $\S_{D}: L^2(\partial D) \rightarrow H^1(\p D)$, defined by
\begin{equation*}
	\S_D[\phi](x) := - \int_{\partial D} \frac{1}{4\pi|x-y|} \phi(y) \dx \sigma(y), \quad x \in \R^3. 
\end{equation*}
It is well-known that, in three dimensions, $\S_D$ is invertible (see, for instance, \cite{ammari2018mathematical}). Here, $H^1(\partial D)$ is the usual Sobolev space of
square-integrable functions on $\partial D$ whose weak derivative is square integrable.

\begin{defn}[capacitance matrix]
We define the basis functions $\psi_i$ and the capacitance coefficients $C_{ij}$ as
\begin{equation} \label{capdef} \psi_i = \left(\S_D\right)^{-1}[\chi_{\p D_i}], \qquad C_{ij}= -\int_{\p D_i} \psi_j   \dx \sigma,\end{equation}
for $i,j=1,...,N$, where $\chi_{\partial D_i}$ is the characteristic function of $\partial D_i$.  The {\em capacitance matrix} $C$ is defined as the matrix $C = \left(C_{ij}\right)$. \end{defn}

The following result is from \cite{ammari2020time}.
\begin{thm}
	Assume that the material parameters are given by \eqref{eq:resonatormod_finite}. Then, as $\delta \to 0$, the quasifrequencies $\omega \in Y^*_t$ to the wave equation \eqref{eq:wave} in the subwavelength regime are, to leading order, given by the quasifrequencies of the system of ordinary differential equations
	\begin{equation}\label{eq:C_ODE_finite}
		\sum_{j=1}^N C_{ij} c_j(t) = \frac{|D_i|\rho_r}{\delta\kappa_r}\frac{1}{\rho_i(t)}\frac{\dx}{\dx t}\left(\frac{1}{\kappa_i(t)}\frac{\dx \rho_ic_i}{\dx t}\right),
	\end{equation}
	for $i=1,...,N$.
\end{thm}

We remark that the left-hand side of \eqref{eq:C_ODE_finite} is specified by the entries of the matrix-vector product $C c(t)$. In fact, we can rewrite the leading order of \eqref{eq:C_ODE_finite} into the following system of Hill equations
\begin{equation}	\label{eq:hill}
\Psi''(t)+ M(t)\Psi(t)=0,
\end{equation}
where $\Psi$ is the vector defined as 
$$\Psi = \left(\frac{\rho_i(t)}{\sqrt{\kappa_i(t)}}c_i(t)\right)_{i=1}^N$$
and $M$ is the matrix defined as 
\begin{equation}\label{defM}
M(t) = \frac{\delta \kappa_r}{\rho_r}W_1(t)C W_2(t) + W_3(t)
\end{equation}
with $W_1, W_2$ and $W_3$ being the diagonal matrices with diagonal entries 
$$\left(W_1\right)_{ii} = \frac{\sqrt{\kappa_i}\rho_i}{|D_i|}, \qquad \left(W_2\right)_{ii} =\frac{\sqrt{\kappa_i}}{\rho_i}, \qquad \left(W_3\right)_{ii} = \frac{\sqrt{\kappa_i}}{2}\frac{\dx }{\dx t}\frac{\kappa_i'}{\kappa_i^{3/2}},$$
for $i=1,...,N$.

Let us remark that multiplying the capacitance matrix $C$ from the left with $W_1$ is equivalent to magnifying each row of the capacitance matrix by a factor constituted by the corresponding element of $W_1$ (as the latter is diagonal). An analogous effect is obtained from the right multiplication with $W_1$, which this time enhances the columns of $C$. As a result, $\rho$ cancels in the diagonal components of $M$:

\begin{equation}\label{diagonal}
    M_{ii} = \frac{\delta \kappa_r \kappa_i}{\rho_r |D_i|} C_{ii} + (W_3)_{ii}.
\end{equation}
This will help us to provide a physical interpretation of the Hamiltonian form in the next section.

\begin{rmk}
In the case of an infinite periodic lattice, we could have used Floquet-Bloch theory (see \cite{ammari2020time}) to reduce ourselves to studying a single cell with a finite number of resonators. In this way, in this paper the infinite lattice case could have been dealt with in an analogous way as for the case of finite number of resonators. 
\end{rmk}

\section{Hamiltonian form}
In the remainder of this paper, we will study the capacitance matrix formulation (\ref{eq:C_ODE_finite}), which can be viewed as a discrete approximation of  the wave equation (\ref{eq:wave}) for computing subwavelength quasifrequencies.

We first focus on transforming the initial system of Hill equations \eqref{eq:hill} into a first-order homogeneous system. Our first approach is to write
\begin{equation}
    i
    \begin{pmatrix}
    \Psi(t) \\
    \Psi ' (t)
    \end{pmatrix}' =
    i\begin{pmatrix}
    0 & \mathbb{Id}_N \\
    -M(t) & 0 
    \end{pmatrix} \begin{pmatrix}
    \Psi(t) \\
    \Psi ' (t)
    \end{pmatrix},
\end{equation}
where $'$ denotes differentiation with respect to time and $\mathbb{Id}_N$ denotes the identity matrix of dimension $N$. 

We will henceforth refer to the matrix of the coefficients of this system as a Hamiltonian. For our purposes, we would like the Hamiltonian to be Hermitian. The main motivation for this, as we will see in the next section, is to show the differences and the similarities between the capacitance matrix formulation \eqref{eq:C_ODE_finite} and the tight-binding approximation in condensed matter theory. Furthermore, this Hamiltonian formulation may also act as starting structure to establish an effective medium theory for large systems of subwavelength modulated resonators \cite{ammari2017effective}. Namely, it will consist in obtaining an estimate for equations \eqref{eq:wave} and \eqref{eq:hill} as the number of the resonators tends to infinity. In all of these applications, it will also be important to us to determine a procedure to obtain the capacitance matrix from the Hamiltonian and vice versa. The main motivation for this is to keep the physical meaning of the components of the two matrices. While for the capacitance matrix formulation the $i,j$-th component represents the magnitude of the interaction between the $i$-th and $j$-th resonators ($i,j =1,..,N$), the physical meaning is less clear to read from the Hamiltonian. One could try to interpret the $i,j$-th component of the latter in a ``kinematic way'' by viewing it as the contribution of the $j$-th ``position'' (or ``velocity'' depending on whether $j$ is greater than $N$ or not) of the value in the resonators to the $i$-th ``position'' or ``velocity''. A more proficuous approach, as we will see, will be to identify a way to go back to the capacitance matrix formulation.
\subsection{Setting}

We are therefore looking for a $\tau$-periodic linear invertible transformation $T(t)$ such that the transformed Hamiltonian matrix is Hermitian. The transformed equation reads 
\begin{equation} \label{eq:Hamiltonian}
     i \Phi' (t) = H (t) \Phi (t),   
\end{equation}
where
\begin{equation}\label{transformed solution to Hill equation}
    \Phi (t) := T(t)  
    \begin{pmatrix}
    \Psi(t) \\
    \Psi ' (t)
    \end{pmatrix}
\end{equation}
and 
\begin{equation} \label{MatrixDefH}
    H(t) := i \left ( T'(t) + T(t) \begin{pmatrix}
    \mathbf{0} & \mathbb{Id}_N \\
    -M(t) &\mathbf{0} 
    \end{pmatrix}\right) T^{-1}(t).
\end{equation}
As previously emphasised, the main condition our Hamiltonian $H$  should satisfy is Hermitianity. We therefore impose the Hermitianity conditions component-wise and obtain a system of equations for the components of $T$. Hence,

$$
    H_{j,k} = \overline{H_{k,j}} \ \ \ \text{for }1\leq j \leq k \leq 2N,  
$$

reads 
\begin{align}\label{General Component-wise}
    \sum_{p=1}^{N}[(T'_{j,p} - \sum_{l=1}^{N} T_{j, N+l} M_{l,p}) T^{-1}_{p,k}] + \sum_{p=N+1}^{2N}(T'_{j,p} + T_{j,p-N})T^{-1}_{p,k}  \nonumber \\
    = -\sum_{p=1}^{N}[(\bar{T}'_{k,p} - \sum_{l=1}^{N} \bar{T}_{k, N+l} M_{l,p}) \bar{T}^{-1}_{p,j}] + \sum_{p=N+1}^{2N}(\bar{T}'_{k,p} + \bar{T}_{k,p-N}) \bar{T}^{-1}_{p,j}, 
\end{align}
for $1\leq j\leq k \leq 2N$, where $N$ denotes the number of resonators.

Notice that we obtain in this way $N(2N+1)$ complex equations. That is because if we consider the equations arising by imposing the lower triangular components of $H$ to be equal to the conjugate of its corresponding upper triangular component, and vice versa, we would obtain the same equation. Here, without loss of generality, we have imposed the condition on the upper triangular components.

By supposing that $T$ is invertible, we will always suppose that $T^{-1}$ is a function of the components of $T$. The purpose of many of the calculations in this section will be to make this relationship explicit. At first, hoping that the equations can be solved in their original form, we will use Cramer's rule to achieve this. Upon realising that the system is too complicated to be efficiently solved in this way (even in the case of a single resonator!), we will employ asymptotic analysis tools together with the Neumann series of a linear operator (which will be truncated at first-order).

We will suppose the transformation $T$ to be of the form
$$
    T(t) = T^0 + \varepsilon T^1(t), \ \ \ t \in [0,\infty), \ \varepsilon >0, 
$$
where $\varepsilon$ denotes the amplitude of the modulation and the time-dependent part $T^1$ will contribute asymptotically less (because of the linearly asymptotic factor $\varepsilon$) than the $T^0$-part. The latter will correspond to the transformation for the case where the parameter modulation is constant in time. Analogously, the modulation of $\rho_i$ and $\kappa_i$, $1 \leq i \leq N$ will be supposed asymptotically small so that the coefficient matrix $M$ defined in \eqref{defM} is of the form
$$
    M(t) = M^0 + \varepsilon M^1(t),
$$
where $M^0$ is the generalised capacitance matrix given by \cite{ammari2021functional}

\begin{equation}\label{defM^0}
    M^0_{ij} := \delta \frac{\kappa_r}{\rho_r |D_i|} C_{ij}.
\end{equation}
Here, $C$ stands for the capacitance matrix introduced in (\ref{capdef}). As a simplifying assumption, we will suppose that all the resonators have the same size. This way, the generalised capacitance matrix \eqref{defM^0} is symmetric positive definite \cite{ammari2021functional}.

As we will see, the Hamiltonian will take on the form
$$
    H(t) = H^0 + \varepsilon H^1(t), \ \ \ t \in [0,\infty), \ \varepsilon >0, 
$$
where the time-dependent part $H^1$ will contribute asymptotically less (through the 
first-order approximation in $\varepsilon$) than the $H^0$-part. $H^0$ will correspond to the Hamiltonian for the case where the material parameters are  constant in time. 

We would like to bring the reader's attention to the fact that the system is under-determined: the number of unknowns, which are given by the components of $T$, is $4N^2$. This is greater than the number of equations, which is only $N(2N+1)$. This suggests that we will need to impose extra conditions to obtain a solution. In the static case, a reasonable extra condition that could be imposed is to suppose the Hamiltonian $H$ to be of the form
\begin{equation}\label{eq:form}
H=\begin{pmatrix}
\mathbf{0} & \mathbf{A} \\
\mathbf{A} & \mathbf{0}
\end{pmatrix}, 
\end{equation}
where  $\mathbf{A}$ is a $N \times N$ symmetric non-negative matrix. The capacitance matrix can be recovered from it by squaring $\mathbf{A}$: we will see that $\mathbf{A}$ is the square root of the generalised capacitance matrix in the static case (subsection below). Therefore, we are able to give a physical meaning to the components of our Hamiltonian. 

Apart from retaining a physical meaning and a connection with the capacitance matrix, we will not be interested in a precise solution. \emph{A-priori}, our search may look somewhat arbitrary, as we will only impose Hermitianity conditions and ``guess'' extra conditions that enable our system to be solvable. Later, we will see that our choice actually also satisfies the other conditions we sought. The extra constraints we will impose will only be reducing the number of variables that come into play (by setting certain components of the transformation $T$ to be 0), instead of imposing extra conditions on $H$.

\subsection{Static case}
In this section, we suppose that the material parameters inside the resonators are fixed. In this case, energy should be conserved, and we seek a Hamiltonian of the form \eqref{eq:form}.

We will denote by $T_0$ and $H_0$ the transformation matrix and the Hamiltonian in the static case respectively. This case is much simpler, as the equation \eqref{MatrixDefH} boils down to
\begin{equation}\label{eq:staticsys}
    H^0 := i \left (T^0 \begin{pmatrix}
    \mathbf{0} & \mathbb{Id}_N \\
    -M^0 & \mathbf{0}
    \end{pmatrix}\right) \left(T^0 \right)^{-1},
\end{equation}
where we recall the definition  \eqref{defM^0} of $M^0$.  Notice that, since we suppose the material parameters and the size to be the same for each resonator, the matrix $M^0$ is symmetric positive definite. We therefore define $\sqrt{M^0}$ as the unique symmetric positive definite matrix such that $\left (\sqrt{M^0} \right )^2 = M^0$.  Then it can be easily checked that one solution to \eqref{eq:staticsys} is given by the transformation matrix
\begin{equation}\label{static_T}
    T^0 := \begin{pmatrix}
    \sqrt{M^0} &\mathbf{0} \\
    \mathbf{0} & i\mathbb{Id}_N
    \end{pmatrix},
\end{equation}
corresponding to the Hamiltonian
\begin{equation}\label{static_H}
     H^0 := \begin{pmatrix}
    \mathbf{0} & \sqrt{M^0} \\
    \sqrt{M^0} & \mathbf{0}
    \end{pmatrix}.
\end{equation}
Although many other transformations would result in a Hermitian Hamiltonian, this particular choice  yields a Hamiltonian of the desired form \eqref{eq:form}. The remaining off-diagonal blocks can be multiplied together to obtain the generalised capacitance matrix. It follows that we are able to give a meaning to the elements of the Hamiltonian $H^0$: 
\begin{equation}\label{Phys_Interpretation_Ham_static}
    M^0_{ij} = \sum_{k=1}^N H^0_{i,N+k} H^0_{N+k,j}.
\end{equation}
To obtain the equation \eqref{Phys_Interpretation_Ham_static} we have also used that, when the material parameters are constant in time, $W_3 \equiv 0$ and the material parameters are the same for all the resonators. The relation \eqref{Phys_Interpretation_Ham_static} states that the $i$-th row of the non-zero blocks of the Hamiltonian contributes to the $i$-th row of the generalised capacitance matrix, and analogously for the columns. Since the $ij$-th component of the capacitance matrix represents the interaction between the $i$-th and $j$-th resonators, the $i$-th row (or equivalently column by symmetry) accounts for the interactions of the $i$-th resonator with the rest of the system, via the formula \eqref{Phys_Interpretation_Ham_static}.

\subsection{Block form of the Hamiltonian}

As the computations will be involved, it will be useful to cast equations \eqref{General Component-wise} in an alternative form. This form will also be useful to determine the compatibility of the system (and hence prove existence of a solution through the standard Picard-Lindel\"of theorem for existence and uniqueness of solutions  to systems of ODEs in normal form with smooth coefficients) for an arbitrary number of resonators. 

We intend to rewrite the equation  \eqref{MatrixDefH} for $H$ completely in block form: it should be a $2 \times 2$ block matrix, where each block is $N\times N$. In order to do so, we need an explicit formula for the inverse of T.

For this, let us write
$$
T:=\begin{pmatrix}
T_1 & T_2 \\ T_3 & T_4
\end{pmatrix}.
$$

As previously mentioned, we will employ perturbative methods (Taylor series truncated at first-order in $\varepsilon$) and suppose $T(t) =T^0+ \varepsilon T^1$. Now, we will use Neumann series truncated at first-order in $\varepsilon$ 
to obtain
\begin{align}\label{T_inverse_Neumann}
        T^{-1} &= \left (T^0 + \varepsilon T^1 \right )^{-1} 
        = (T^0)^{-1} \left [\mathbb{Id}_N - \left ( - \varepsilon T^1 (T^0)^{-1} \right)  \right ]^{-1}  \nonumber \\
        & =  (T^0)^{-1} \sum_{k=0}^\infty \left ( - \varepsilon T^1 (T^0)^{-1} \right) \approx (T^0)^{-1} - \varepsilon (T^0)^{-1} T^1 (T^0)^{-1}.
\end{align}

Since the matrix for the transformation in the static case $T^0$ is block diagonal (and all the non-zero blocks are invertible), then it is possible to easily compute its inverse through the identities regarding the inverse of a $2\times 2$ block matrix. We obtain
$$
     (T^0)^{-1} = \begin{pmatrix}
     \left( M^0 \right) ^{-\frac{1}{2}} & \mathbf{0} \\
       \mathbf{0}  & -i \mathbb{Id}_N
     \end{pmatrix}.
$$

Plugging this into the equation \eqref{MatrixDefH} for $H$  and noticing that by differentiating $T$ in perturbed form only the first-order term survives, we obtain 

\begin{align}\label{Block_Form_H}
 &   H^0 + \varepsilon H^1(t) \nonumber \\ &
    \approx i \left ( T'(t) + T(t) \begin{pmatrix}
    \mathbf{0} & \mathbb{Id}_N\\
    -M(t) & \mathbf{0} 
    \end{pmatrix}\right) \left ((T^0)^{-1} - \varepsilon (T^0)^{-1} T^1 (t) (T^0)^{-1} \right) \nonumber\\&
    = i \Biggl ( \varepsilon \begin{pmatrix}
(T^1_1)' & (T_2^1)' \\ (T_3^1)' & (T_4^1)'
\end{pmatrix}(t) + \left [ \begin{pmatrix}
    \sqrt{M^0} &\mathbf{0} \\
    \mathbf{0} & i\mathbb{Id}_N
    \end{pmatrix} + \varepsilon \begin{pmatrix}
    T_1^1 & T_2^1 \\ T_3^1 & T_4^1
    \end{pmatrix}(t) \right ] \begin{pmatrix}
    \mathbf{0} & \mathbb{Id}_N\\
    -M(t) & \mathbf{0} 
    \end{pmatrix}  \nonumber\\ &
    \cdot \left (\begin{pmatrix}
     \left( M^0 \right) ^{-\frac{1}{2}} & \mathbf{0} \\
       \mathbf{0}  & -i \mathbb{Id}_N
     \end{pmatrix}- \varepsilon \begin{pmatrix}
     \left( M^0 \right) ^{-\frac{1}{2}} & \mathbf{0} \\
       \mathbf{0}  & -i \mathbb{Id}_N
     \end{pmatrix}\begin{pmatrix}
    T_1^1 & T_2^1 \\ T_3^1 & T_4^1
    \end{pmatrix}(t) \begin{pmatrix}
     \left( M^0 \right) ^{-\frac{1}{2}} & \mathbf{0} \\
       \mathbf{0}  & -i \mathbb{Id}_N
     \end{pmatrix} \right) \nonumber \\&
     = H^0 +\nonumber \\&
     \varepsilon \begin{pmatrix}
     i(T^1_1)' (M^0)^{-\frac{1}{2}} - i T^1_2 \sqrt{M^0} - \sqrt{M^0}T^1_3 (M^0)^{-\frac{1}{2}} & (T^1_2)' + T^1_1 +i \sqrt{M^0} T^1_4 \\
     \nm
     i(T^1_3)' (M^0)^{-\frac{1}{2}} - i T^1_4 \sqrt{M^0} +(M^0)^{-\frac{1}{2}}M^1 - \sqrt{M^0}T^1_1 (M^0)^{-\frac{1}{2}}  &(T^1_4)' + T^1_3 +i \sqrt{M^0} T^1_2 \end{pmatrix}.
\end{align}

We already know that $M^0$ is real. As an additional simplification, we will also suppose the modulated part $M^1 (t)$ to be real for all $t \geq 0$. This means that energy gain/loss only comes from the time-varying part. Otherwise, we should not even expect the Hamiltonian to be Hermitian. 

\subsection{One-resonator case}
This subsection and the following will be devoted to finding a solution to the under-determined system of equations \eqref{General Component-wise} (or its asymptotic counterpart given by the block notation \eqref{Block_Form_H}).  As is clear from the block form \eqref{Block_Form_H}, (recalling that the equations arise by imposing the matrix to be Hermitian), the system may be incompatible due to the inhomogeneity given by $M$. Finding a solution can be challenging, and correct additional assumptions must be made to obtain a compatible system.

We focus on computing the $2\times2$ matrix $T$ in the case of a single modulated resonator. The system \eqref{General Component-wise} of three equations reads
\begin{equation*}
    \begin{cases}
  \ds  \mathcal{R}\left((T'_{1,1} - T_{1,2} M) T^{-1}_{1,1} + (T'_{1,2} + T_{1,1})T^{-1}_{2,1}\right) = 0 &\text{ if } (j,k) = (1,1),\\
    \nm \ds
    \mathcal{R}\left((T'_{2,1} - T_{2,2} M) T^{-1}_{1,2} + (T'_{2,2} + T_{2,1})T^{-1}_{2,2}\right) = 0 &\text{ if } (j,k) = (2,2), \\
    \nm \ds
    (T'_{1,1} - T_{1,2} M) T^{-1}_{1,2} + (T'_{1,2} + T_{1,1})T^{-1}_{2,2} \\ 
    \nm \ds = (-\bar{T}'_{2,1} + \bar{T}_{2,2} M) \bar{T}^{-1}_{1,1} - (\bar{T}'_{2,2} + \bar{T}_{2,1}) \bar{T}^{-1}_{2,1} &\text{ if } (j,k)=(1,2),
    \end{cases}
\end{equation*}
where $\mathcal{R}(z)$ denotes the real part of $z \in \mathbb{C}$. 
We may use Cramer's rule to compute $T^{-1}$ as a function of $T$
\begin{equation*}
    \begin{cases}
    \ds T^{-1}_{1,1} = \frac{T_{2,2}}{\det T}, \\
    \nm \ds
    T^{-1}_{1,2} = -\frac{T_{1,2}}{\det T}, \\
    \nm \ds
    T^{-1}_{2,1} = -\frac{T_{2,1}}{\det T},\\
    \nm \ds
    T^{-1}_{2,2} = \frac{T_{1,1}}{\det T}.
    \end{cases}
\end{equation*}
The system then becomes
\begin{equation*}
    \begin{cases}
   \ds \mathcal{R}\left((T'_{1,1} - T_{1,2} M) \frac{T_{2,2}}{\det T} - (T'_{1,2} + T_{1,1})\frac{T_{2,1}}{\det T}\right) = 0 &\text{ if } (j,k) = (1,1),\\
    \nm \ds
    \mathcal{R}\left((-T'_{2,1} + T_{2,2} M)\frac{T_{1,2}}{\det T} + (T'_{2,2} + T_{2,1})\frac{T_{1,1}}{\det T}\right) = 0 &\text{ if } (j,k) = (2,2), \\
    \nm \ds
    (-T'_{1,1} + T_{1,2} M)\frac{T_{1,2}}{\det T} + (T'_{1,2} + T_{1,1})\frac{T_{1,1}}{\det T} \\[0.5em] \ds = (-\bar{T}'_{2,1} + \bar{T}_{2,2} M) \frac{\bar{T}_{2,2}}{\det \bar{T}} + (\bar{T}'_{2,2} + \bar{T}_{2,1}) \frac{\bar{T}_{2,1}}{\det \bar{T}} &\text{ if } (j,k)=(1,2).
    \end{cases}
\end{equation*}
We first attempt to compute a simplified perturbed version of our equations which will turn this into a system of linear ODEs. We suppose $T = T^0 + \varepsilon \ T^1(t)$, for $\varepsilon>0 $ small. We have shown that the $T^0$ associated to the static case corresponds the the $\varepsilon$-zero-order. Let us focus on the first-order equations we obtain. 

The following expression will be useful:
\begin{align*}
   \ds \frac{1}{\det T} = & \frac{1}{\det T^0} - \frac{\varepsilon}{(\det T^0)^2}(T_{1,1}^0 T^1_{2,2} + T_{2,2}^0 T^1_{1,1} - T_{2,1}^0 T^1_{1,2} - T_{1,2}^0 T^1_{2,1})  + O(\epsilon^2) \\
   \nm \ds
    =& -\frac{i}{\sqrt{M^0}} + \frac{\varepsilon}{M^0}(i \ T_{1,1}^1 + \sqrt{M^0} T_{2,2}^1) + O(\epsilon^2),
\end{align*}   
where in the last equality we have used the known values for $T^0$. 

We obtain the following system:

\begin{equation}\label{1res_complex_system}
    \begin{cases}
  \ds  \mathcal{R} \left(\sqrt{M^0} [- T_{1,1}^1 + M^0 T_{1,2}^1 -i \sqrt{M^0} T^1_{2,1} ]\right)=0 \text{ if } (j,k)=(1,1), \\
  \nm \ds
    -i\left((T_{1,2}^1)' + T^1_{1,1}\right) + \sqrt{M^0} T_{2,2}^1 \\ = \sqrt{M^0} \bar{T}_{2,2}^1 + \frac{1}{\sqrt{M^0}} (\bar{T}_{2,1}^1)' - \frac{i}{\sqrt{M^0}} M^1 + i \bar{T}^1_{1,1}  \text{ if } (j,k)=(1,2), \\
    \nm \ds
    \mathcal{R} \left(\sqrt{M^0} [- T_{1,1}^1 + M^0 T_{1,2}^1 -i \sqrt{M^0} T_{2,1} ]\right)=0 \text{ if } (j,k)=(2,2). 
    \end{cases}
\end{equation}

We now intend to rewrite the whole system in real equations by splitting real and imaginary parts. To do so and simplify notation, we write 
\begin{align*}
  T^1 (t) = \begin{pmatrix}
  a + i b & c + i d \\
  e + i f & g + i h
  \end{pmatrix} (t),
\end{align*}
where $a,\ b, \ c,\ d,\ e,\ f,\ g,\ h$ are real-valued functions of time.

The system then becomes:

\begin{equation}\label{1res_real_system}
\begin{cases}
\ds a'= \sqrt{M^0} f +M^0 c,\\
\nm
\ds
h' = -f-\sqrt{M^0}c ,\\
\nm \ds
c' = 2 \sqrt{M^0} h + e' + \frac{M^1}{\sqrt{M^0}} - 2a, \\
\nm \ds
d' = f'.
\end{cases}
\end{equation}

To obtain a solution, we will suppose the modulation to be zero at $t=0$, hence we will impose the initial conditions to be zero for all the unknowns. This way, from the fourth equation in the system we obtain $d=f$.
Moreover, let us  notice that $b$ and $g$, because of the nature of the initial equations, do not appear in our system. 

Let us also notice once again that the system is under-determined: we will impose extra assumptions on our unknowns for the system to have a unique explicit solution. 
The most natural assumption that does not lead (through the inhomogeneity in the third equation) to an incompatible system is the following. We impose $b\equiv g\equiv 0$, as they do not appear in the system.  Moreover, we impose $f$, and hence $d$, to be equal to $0$. Lastly, we also suppose $e$ to be zero. We then obtain the following form for the modulated part of our transformation:
\begin{equation} \label{T_Upper_triangular_N1}
    T^1= \begin{pmatrix}
    a & c \\
    0 & ih
    \end{pmatrix}
\end{equation}
hence it assumes an upper-triangular structure. Here, the diagonal entries are alternatively purely real and imaginary. By imposing $c \equiv 0$ instead of $e$, we would have obtained a similar system with an analogous solution, and $T^1$ would have been lower triangular.

Although an explicit solution for $T^1$  is in general fairly complicated, we can still plug in the equations \eqref{1res_real_system} we obtained into \eqref{Block_Form_H} and notice that the following cancellations occur 
\begin{align}\label{H_with_solution_N1}
    H^1(t) &= \begin{pmatrix}
    i\frac{a'}{\sqrt{M^0}} - i \sqrt{M^0} c & c' + a - \sqrt{M^0} h \\
    \sqrt{M^0} h + \frac{M^1}{M^0} -a & ih' + i \sqrt{M^0}  c
    \end{pmatrix}  \nonumber \\
    &= \begin{pmatrix}
    0 & \sqrt{M^0} h + \frac{M^1}{M^0} -a \\
    \sqrt{M^0} h + \frac{M^1}{M^0} -a & 0
    \end{pmatrix}.
\end{align}

Remarkably, we see that our chosen transformation $T^1$ given by \eqref{T_Upper_triangular_N1} has the effect of yielding a modulated Hamiltonian which is also real block off-diagonal, as in \eqref{eq:form}. The fact that the matrix is purely real is in line with the physical fact that we are supposing the only energy gain/loss to come from the resonator's modulation. To retain physical meaning though, we can no longer look at the capacitance matrix, as with the current definition it only accounts for the geometry of the system \cite{ammari2021functional}. We would have to seek a connection with the system  \eqref{eq:hill} of Hill equations. This is given by the transformation $T$. Indeed, from \eqref{MatrixDefH} we may invert $T$ and compute the initial solution of the Hill equation from the solution of the Hamiltonian system using the relation \eqref{transformed solution to Hill equation}. Moreover, by inverting \eqref{MatrixDefH}, we arrive at
$$
  \begin{pmatrix}
        \mathbf{0} & \mathbb{Id}_N \\
        -M(t) &\mathbf{0} 
        \end{pmatrix}= -i T^{-1}(t) \left (H(t)T(t) - T'(t) \right), 
$$
from which we can relate $H$ to $M$. Notice that this procedure can be done also in the static case and yields  exactly the same result. 
\begin{rmk}
Note that if we assume that $M^1(t)$ has a simple form, say for instance a cosine function, then it is possible  to obtain an explicit solution by plugging in the system \eqref{1res_real_system} in Wolfram Mathematica. Nevertheless, the obtained expressions for $a(t), c(t)$, and $h(t)$ are rather complicated.  
\end{rmk}

\subsection{Case of an arbitrary number of resonators}
Since for the one-resonator case an explicit solution was already hard to obtain in the perturbed case, the main focus of this subsection will be to show at least the existence of a solution to the system \eqref{General Component-wise}. Following the previous case $N=1$, we will make assumptions on the transformation $T^1$ which yields a compatible system with a unique solution. Moreover, we will interpret qualitatively the Hamiltonian we obtain, and compare it with what we were expecting to obtain based on a physical interpretation. 

The one-resonator case was mainly studied to gain some intuition on the structure of the more general setting. Indeed, we will see that an extension of that idea will prove to be fruitful in obtaining suitable assumptions on our transformation. Throughout the subsection, we will suppose $T^0$ to be of the form we determined in \eqref{static_T}.

The system of equations \eqref{General Component-wise} can be divided, using  the block notation \eqref{Block_Form_H} for $H$, into three subsystems, obtained by verifying what conditions the Hermitianity of $H$ imposes on the two diagonal blocks and on the off-diagonal ones, respectively. In particular, imposing $H$ to be Hermitian is equivalent to imposing the diagonal blocks of \eqref{Block_Form_H} to be Hermitian, and the $(j,$ $k)$ component of the upper off-diagonal block to be equal to the conjugate of the $(k,$ $j)$ component of the lower off-diagonal block, for $1 \leq j,k \leq N$. More explicitly,  we have

\begin{enumerate}
    \item Hermitianity of $H$ implies Hermitianity of the diagonal blocks. Imposing this on the first diagonal block in \eqref{Block_Form_H}, we obtain:
    \begin{align} \label{first_block_system}
        \sum_{l=1}^N \left \{i \left (T^1_{jl} \right)' \left (M^0 \right)^{-\frac{1}{2}}_{lk} - i \left (T^1_{j,l+N} \right) \sqrt{M^0}_{lk} - \sqrt{M^0}_{jl} \sum_{r=1}^N\left( T^1_{l+N, r} \right) \left (M^0 \right)^{-\frac{1}{2}}_{rk} \right\}  \nonumber \\
        = \sum_{l=1}^N \left \{- i \left (\bar{T}^1_{kl} \right)' \left (M^0 \right)^{-\frac{1}{2}}_{lj} + i \left (\bar{T}^1_{k,l+N} \right) \sqrt{M^0}_{lj} - \sqrt{M^0}_{kl} \sum_{r=1}^N\left( \bar{T}^1_{l+N, r} \right) \left (M^0 \right)^{-\frac{1}{2}}_{rj} \right\},
    \end{align}
    where we let $1 \leq j \leq k \leq N$, and hence the obtained system comprises of $N(N+1)$ real equations (by equating real and complex parts of our equation separately). \\
    \item We now impose the other diagonal block to be Hermitian to obtain:
    \begin{align}\label{second_block_system}
        \left ( T^1_{j+N,k+N} \right)' +\left ( T^1_{j+N, k}\right) + i \sum_{l=1}^N \sqrt{M^0}_{jl} \left ( T^1_{l,k+N} \right )  \nonumber \\
        =  \left ( \bar{T}^1_{k+N,j+N} \right)' +\left ( \bar{T}^1_{k+N, j}\right) - i \sum_{l=1}^N \sqrt{M^0}_{kl} \left ( \bar{T}^1_{l,j+N} \right ),
    \end{align}
    where again $1 \leq j \leq k  \leq N$, hence the obtained system comprises of $\frac{N(N+1)}{2}$ complex equations. \\
    \item Lastly, let us impose the upper off-diagonal block to be equal to the transpose of the complex adjoint of the lower off-diagonal block:
    \begin{align}\label{off-diagonal_block_system}
        & \left ( T^1_{j,k+N} \right)' +\left ( T^1_{j, k}\right) + i \sum_{l=1}^N \sqrt{M^0}_{jl} \left ( T^1_{l+N,k+N} \right ) =  \sum_{l=1}^N \biggl \{ - i \left (\bar{T}^1_{k+N,l} \right)' \left (M^0 \right)^{-\frac{1}{2}}_{lj}  + \nonumber \\ & + i \left (\bar{T}^1_{k+N,l+N} \right) \sqrt{M^0}_{lj}+ \left (M^0 \right)^{-\frac{1}{2}}_{kl} M^1_{lj} - \sqrt{M^0}_{kl} \sum_{r=1}^N\left( \bar{T}^1_{l, r} \right) \left (M^0 \right)^{-\frac{1}{2}}_{rj} \biggr\},
    \end{align}
    where this time instead we have $1 \leq j,k \leq N$. Therefore, the system has $2N^2$ real equations.
\end{enumerate}
One can easily see that the system is under-determined as there are $4N^2$ components of $T^1$ while the complex equations are $N(2N +1)$.
To seek the conditions which make the system determined and compatible, we generalise our ``guess'' for the one-resonator case in the following way: we suppose $T$ to be upper triangular, with the first $N$ diagonal elements being purely real and the last $N$ being purely imaginary. This structure for the diagonal elements is justified by the structure of the equations arising by imposing the diagonal elements of the Hamiltonian to be real. Indeed, from systems \eqref{first_block_system} and \eqref{second_block_system} respectively we obtain 
\begin{equation}\label{first diagonal elements}
    \mathcal{R}\left( \sum_{l=1}^N \left \{i \left (T^1_{jl} \right)' \left (M^0 \right)^{-\frac{1}{2}}_{lj} - i \left (T^1_{j,l+N} \right) \sqrt{M^0}_{lj} - \sqrt{M^0}_{jl} \sum_{r=1}^N\left( T^1_{l+N, r} \right) \left (M^0 \right)^{-\frac{1}{2}}_{rj} \right\} \right) = 0, 
\end{equation}
where $1\leq j \leq N$, and 
\begin{equation}\label{second_diagonal_elements}
    \mathcal{R} \left ( \left ( T^1_{j+N,k+N} \right)' +\left ( T^1_{j+N, k}\right) + i \sum_{l=1}^N \sqrt{M^0}_{jl} \left ( T^1_{l,k+N} \right )\right) =0,
\end{equation}
for $1\leq j \leq N$. 

These equations have the same structure as the first and last ones of the one-resonator system \eqref{1res_complex_system}. This peculiarity was what determined the system in the one resonator case to be independent of $b,\ g$, which stood for the imaginary part of the first diagonal component and the real part of the second diagonal component, respectively. In an analogous way, it is possible to check that the same thing happens in the general $N$-resonator case from the subsystems \eqref{first diagonal elements} (where we see the imaginary parts of the diagonal entries of $T$ cancelling) and \eqref{second_diagonal_elements} (this time the real part will cancel). Hence, the assumption on the diagonal elements of $T$  is justified.

With this choice of $T$, the systems \eqref{first_block_system}, \eqref{second_block_system} and \eqref{off-diagonal_block_system} become respectively
  \begin{align}\label{new_first_block_system}
        &\sum_{l=j}^N \left [i \left (T^1_{jl} \right)' \left (M^0 \right)^{-\frac{1}{2}}_{lk} \right] - i \sum_{l=1}^N \left [\left (T^1_{j,l+N} \right) \sqrt{M^0}_{lk}\right ]  \nonumber \\
       & = -i \sum_{l=k}^N \left [ \left (\bar{T}^1_{kl} \right)' \left (M^0 \right)^{-\frac{1}{2}}_{lj} \right ] + i \sum_{l=1}^N \left [\left (\bar{T}^1_{k,l+N} \right) \sqrt{M^0}_{lj} \right ],
        \\\label{new_second_block_system}
       & \delta_{j\leq k}  \left ( T^1_{j+N,k+N} \right)' + i \sum_{l=1}^N \sqrt{M^0}_{jl} \left ( T^1_{l,k+N} \right ) \nonumber \\
       & =  \delta_{k \leq j}\left ( \bar{T}^1_{k+N,j+N} \right)' - i \sum_{l=1}^N \sqrt{M^0}_{kl} \left ( \bar{T}^1_{l,j+N} \right ),
\\
\label{new_off-diagonal_block_system}
     & \left ( T^1_{j,k+N} \right)' + \delta_{j < k}\left ( T^1_{j, k}\right) + i \sum_{l=1}^N \sqrt{M^0}_{jl} \left ( T^1_{l+N,k+N} \right )  \nonumber \\
        & = i \sum_{l=k}^N \left [\left (\bar{T}^1_{k+N,l+N} \right) \sqrt{M^0}_{lj} \right] + \sum_{l=1}^N \left [\left (M^0 \right)^{-\frac{1}{2}}_{kl} M^1_{lj} - \sqrt{M^0}_{kl} \sum_{r=l}^N\left( \bar{T}^1_{l, r} \right) \left (M^0 \right)^{-\frac{1}{2}}_{rj} \right].
    \end{align}
\newline  We now suppose $T^1_{jk} (t) := a_{jk}(t) + i b_{jk}(t)$, where $a_{jk},\ b_{jk}$ are real valued functions, for $1 \leq j < k \leq 2N$. For the diagonal entries, we have $T^1_{jj}(t) := a_{jj} (t)$, $1\leq j \leq N$ and $T^1_{jj}(t) := ib_{jj} (t)$, $N+1\leq j \leq 2N$. Let us point out that we did not define either $a_{jj}$ for $j>N$ nor $b_{jj}$ for $j\leq N$. This is in line with our previous remark regarding our choice for the structure of the diagonal entries of $T$. As further confirmation, one can verify that we would obtain the identity $0=0$ in the equations where those variables should appear.

Moreover, without loss of generality we will suppose $M^0$ (and hence $\sqrt{M^0}$ and $\left ( M^0 \right)^{-\frac{1}{2}}$) to be diagonal. By separating real and imaginary parts of the previous system, we then obtain

\begin{align}\label{First_projected general N}
    &\begin{cases}
            a_{jk}' \left (M^0 \right)^{-\frac{1}{2}}_{kk} - a_{j,k+N} \sqrt{M^0}_{kk} = - \delta_{j=k} \ a'_{kj}\left (M^0 \right)^{-\frac{1}{2}}_{jj} + a_{k,j+N} \sqrt{M^0}_{jj}, \\
            \nm
            \delta_{j < k}\  b_{jk}' \left (M^0 \right)^{-\frac{1}{2}}_{kk} - b_{j,k+N} \sqrt{M^0}_{kk} =  - b_{k,j+N} \sqrt{M^0}_{jj},
    \end{cases}
        \\\label{second_projected general N}
    &\begin{cases}
 \delta_{j < k}\ \ a_{j+N,k+N}' - \sqrt{M^0}_{jj} b_{j,k+N} = - \sqrt{M^0}_{kk} b_{k,j+N},\\
 \nm
 b_{j+N,k+N}' + \sqrt{M^0}_{jj} a_{j,k+N}  =  -\delta_{j= k}\ \ b_{k+N,j+N}' - \sqrt{M^0}_{kk} a_{k,j+N},
       \end{cases}
\\
\label{projected_off-diagonal_block_system}
    &\begin{cases}
     (a_{j,k+N})' + \delta_{j \leq  k} \ a_{j k}  - \delta_{j \leq k}\sqrt{M^0}_{jj} b_{j+N,k+N} \\ 
     = \delta_{k \leq j} \ b_{k+N,j+N} \sqrt{M^0}_{jj} +  \left (M^0 \right)^{-\frac{1}{2}}_{kk} M^1_{kj} - \delta_{k \leq j} \sqrt{M^0}_{kk}  a_{kj}  \left (M^0 \right)^{-\frac{1}{2}}_{jj}, \\ \\
     \nm
     b_{j,k+N}' + \delta_{j < k} \  b_{j k}  + \delta_{j < k} \sqrt{M^0}_{jj} a_{j+N,k+N}  \\ 
     \nm
     =\delta_{k < j} \ a_{k+N,j+N} \sqrt{M^0}_{jj} +  \left (M^0 \right)^{-\frac{1}{2}}_{kk} M^1_{kj} + \delta_{k <j} \sqrt{M^0}_{kk}  b_{kj}  \left (M^0 \right)^{-\frac{1}{2}}_{jj},
        \end{cases}
\end{align}
where we let $1 \leq j \leq k \leq N$ in the first two subsystems and $1 \leq j,k \leq N$ in the last one. We couple this with the initial condition $\mathbf{a}(0) = \mathbf{b}(0)=0$, where $\mathbf{a}$ and $\mathbf{b}$ stand for the vector valued functions of time whose components are respectively  $a_{jk}$ and $b_{jk}$  for appropriate $j$, $k$. 

Let us remark that the system is now determined: the number of equations coincides with the number of unknowns not set to zero. The only thing left to verify is that the system is compatible (let us recall that the incompatibility may arise only because of the inhomogeneity in the third subsystem). For this purpose, let us look at the three systems \eqref{First_projected general N}, \eqref{second_projected general N} and  \eqref{projected_off-diagonal_block_system} in matrix form,
$$
B\begin{pmatrix}
	\mathbf{a}'(t) \\
	\mathbf{b}'(t)
\end{pmatrix} = A {\begin{pmatrix}
		\mathbf{a}(t) \\
		\mathbf{b}(t)
\end{pmatrix}} + \mathbf{c}(t),
$$
where $B$ is the $4N^2 \times 4N^2$ matrix of the coefficients of the derivatives, $\mathbf{a} \in \mathbb{R}^{2N^2}$ and $\mathbf{b} \in \mathbb{R}^{2N^2}$ are the vector of the unknowns $(a_{jk})$ and $(b_{ij})$ respectively, and $A$ is the $4N^2 \times 4N^2$ matrix of the coefficients of $\mathbf{a}$ and $\mathbf{b}$. Lastly, $\mathbf{c}$ is the vector of the inhomogeneity whose components would be linear combinations of the components of $M^1$. If $B$ would be invertible, since $\mathbf{c}$ is smooth, we could apply standard existence and uniqueness theorems for first-order systems of ODEs to obtain a solution.

To see that $B$ is invertible, it suffices to notice that our assumptions actually turn each row into a vector with only one entry different from zero, and this entry is different for every row. Indeed, the systems \eqref{First_projected general N}, \eqref{second_projected general N} and \eqref{projected_off-diagonal_block_system} can each be divided in two sub-systems each of which contains either only derivatives of $\mathbf{a}$ or of $\mathbf{b}$. For the first sub-system of \eqref{First_projected general N}, the only coefficients different from zero are the ones of $a_{jk}'$, $1\leq j\leq k \leq N$. For the second one of \eqref{First_projected general N}, the same happens for $b_{jk}'$, $j \not = k$. Analogously, systems \eqref{second_projected general N} and \eqref{projected_off-diagonal_block_system} do the same for the remaining derivatives of the unknowns of the system.

Finally, let us remark that the structure of the Hamiltonian \emph{a-priori} will not have the peculiar structure of the one-resonator case. In general, it will not be purely real nor off-diagonal. However, one can check by plugging in \eqref{Block_Form_H} the equations with $j=k$ in systems \eqref{First_projected general N} and \eqref{second_projected general N} that the diagonal entries of $H$ are $0$. This is in line with the symmetry in this Hamiltonian's eigenvalues: it has both the positive and the negative roots of the capacitance matrix (at leading-order). Moreover, in the same way by using the relations found in the system \eqref{projected_off-diagonal_block_system} that the diagonal entries of the off-diagonal blocks of $H$ are purely real.

\begin{rmk}
    We would like to point out that since $M(t)$ is periodic we can apply the Floquet theorem (see, for instance, \cite[p. 92]{teschl2012ordinary}). We will need to transform our inhomogeneous system into a homogeneous one, which can be done in a standard way by making the system larger. Hence, from the periodicity of the modulation (made explicit in $M$), we know that both the transformation $T$ that solves \eqref{General Component-wise} and the solution of the transformed system \eqref{eq:Hamiltonian} will be quasiperiodic. Since also the solution of the original system \eqref{eq:hill} will be quasiperiodic, its quasifrequencies can be computed from the ones of $T$ and of the system \eqref{eq:Hamiltonian}. This is important because the quasifrequencies of solutions to \eqref{eq:hill} represent the subwavelength quasifrequencies of the system of $N$ modulated subwavelength resonators (see \cite{ammari2020time}).
\end{rmk}

\section{Tight-binding approximation}
In the field of condensed matter theory in some settings a useful assumption is the tight-binding approximation. It consists in assuming that the major interaction of the system comes from self-interaction. The approximating modes are constructed by taking linear combinations of the modes obtained by only considering the self-interaction of a single atom at the time. This can be proven to be a good approximation for the Schr\"odinger equation when the electrons are tightly bonded to the atom's nuclei in terms of eigenvalues (see, for instance, \cite{osti_377103,https://doi.org/10.1002/cpa.21735}). 

Our goal is to compare this to the capacitance matrix approximation in the theory of wave propagation in metamaterials which also consists of a discretisation of a continuous problem. In this section, we will point out similarities and differences between the two models. To enhance this comparison we will work in the Hamiltonian setting discussed in section 2. We will mainly work in the unmodulated case. As we have a clear representation and connection with the capacitance matrix, this will help us check whether our guesses for the approximants are good or not. This type of inference can then be extended to the modulated case: although we know an explicit closed form for the Hamiltonian for a fixed arbitrary number of resonators, such solution can arguably be too complicated to handle. Therefore, a tight-binding-type approximation could help in yielding a less cumbersome  zero-order approximation for the modulated Hamiltonian starting from the one resonator case. For higher-orders, the involved computations cannot be avoided in general.

The main fundamental difference between the tight-binding regime and the capacitance formulation, aside from the fact that they were conceived to model different problems, is the following. As proven in \cite{ammari2021functional}, the subwavelength modes need to be almost constant inside the resonators. As these modes will not in general be constant elsewhere, taking linear combinations of modes would contradict the almost-constant nature of the true modes in the resonators. For a visual intuition of this fact, we refer to Figures \ref{Modes of two different single resonators Fig} and \ref{fig:both}. 

\begin{figure}[h]
\centering
\input{Single-resonator_modes.tikz}
\caption{Modes of two different single resonators. The grey shaded areas represent the two resonators, and different colours have been used for the two different modes relative to each resonator separately. Observe that each mode is approximately constant inside corresponding resonator.}
\label{Modes of two different single resonators Fig}
\end{figure}
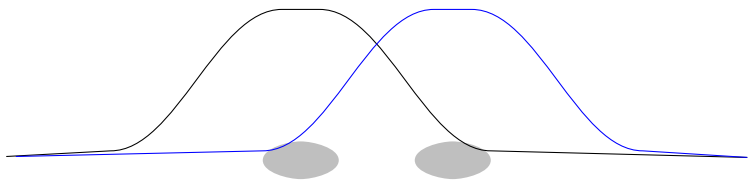

\begin{figure}[h]
	\begin{subfigure}[b]{0.48\linewidth}
\centering
\input{Two-resonator_mode.tikz}
\caption{Example of a mode for two coupled resonators. Observe that the mode is approximately constant inside both resonators.}
\label{Mode when both resonators are considered. Fig}
\end{subfigure}\hfill
	\begin{subfigure}[b]{0.48\linewidth}
\centering
\input{Linear_combination.tikz}
\caption{Linear combination of the modes of the single resonators. \newline}
\label{Linear combination of modes of single resonators. Fig}
\end{subfigure}
\caption{Comparison between a true mode (a) and a tight-binding-type approximant (b). Unless the resonators are asymptotically small, the linear combination in (b) is not constant inside the resonators and is not a suitable approximant of the mode sketched in (a).} \label{fig:both}
\end{figure}
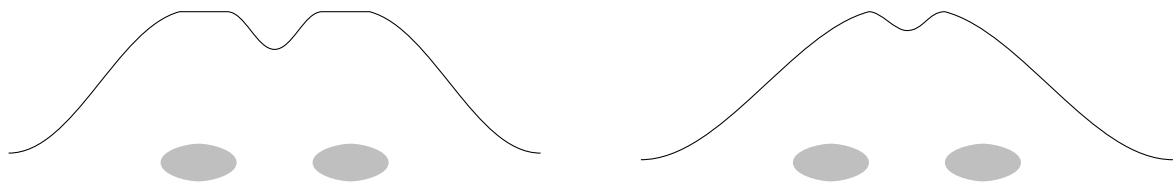

Therefore, for a tight-binding approximation to work in our setting, we will need to assume the size of the resonators to be asymptotically small compared to the distance between them. This regime is called the {\em dilute regime} \cite{ammari2021functional,ammari2017effective}. 

It was shown in \cite{ammari2021functional} that the capacitance matrix has the following form in the dilute regime.
\begin{thm}[capacitance matrix in the dilute regime]
Consider a dilute system of $N$ identical subwavelength resonators, given by
$$
    D= \bigcup_{j=1}^N (B + \eta^{-1} z_j),
$$
where $0<\eta \ll 1$, $B$ is some fixed domain and $\eta^{-1} z_j$ represents the position of each resonator, for fixed $z_j$. In the limit as $\eta \rightarrow 0$, the asymptotic behaviour of the corresponding capacitance matrix is given by
\begin{equation}\label{capacitance matrix in dilute regime}
   \ds C_{ij} = \begin{cases}
    \capacity_B + O(\eta^2), & i=j, \\
    \nm \ds
    - \frac{\eta (\capacity_B)^2}{4 \pi |z_i - z_j|} + O(\eta^2), & i \not = j,
    \end{cases}
\end{equation}
where $\ds \capacity_B := - \int_{\partial B} \mathcal{S}_B^{-1} \left [ \chi_{\partial B} \right ] \dx \sigma$.
\end{thm}

It can be noted that, up to first-order, the diagonal entries of the capacitance matrix in the dilute regime only depend on one resonator and represent the self-interaction of a particular resonator. The off-diagonal entries represent the mutual interactions between two resonators and, in the dilute regime, the entry $C_{ij}$ depends solely on the $i$-{th} and the $j$-{th} resonator. In some sense, the capacitance matrix in dilute regime is a good tight-binding-type approximation as the mutual interactions are represented by lower-order terms.
\begin{rmk}
In \cite{cassier21}, a different setting is studied in which the contrast parameter $\delta \rightarrow \infty$. In this case, which is no longer a subwavelength regime, a tight-binding approximation is expected to hold. Indeed, as the modes are almost constant in the bulk, by taking linear combinations of modes with respect to different resonators we would not encounter the same issue as in our case.
\end{rmk}

We seek a better connection to the classical tight-binding approximation by recasting our problem in Hamiltonian form. The requirements for an approximant to be a good competitor are the following:
\begin{itemize}
    \item Its eigenvalues are equal to the ones of the full Hamiltonian $H$, up to first-order in the dilute parameter $\eta$. Here, the dilute assumption is being used to provide an approximation of the eigenvalues of the Hamiltonian.\\
    \item  The diagonal entries of the tight-binding approximant Hamiltonian should represent the magnitude of the modes related to self-interaction of single resonators. The off-diagonal entries of the approximant should represent the mutual interactions between the resonators. In accordance with the capacitance matrix in dilute regime, these entries should be of higher-order in the dilute regime parameter $\eta$ compared to the diagonal entries. Furthermore, we will see that the off-diagonal terms decay too slowly to be omitted. Therefore, there will also be long-range interactions.\\
    \item We will see that the eigenvalues of the Hamiltonian matrix will be the square roots with sign of the eigenvalues of the capacitance matrix. These would then correspond to modes solving the nonlinear eigenvalue problem in the static case corresponding to \eqref{eq:freq}, see \cite{ammari2021functional}. Each positive eigenvalue of the Hamiltonian would have its negative symmetric counterpart corresponding to the negative part of the same root. The negative modes are artificial and do not bear any physical meaning. Therefore, our approximant should not include any interaction between the positive modes and the negative non-physical ones.\\
\end{itemize}

The approximant will be represented as the sum of two matrices. We let $\widetilde{H}$ to be the zero-th order approximation which will denote the self-interaction of the resonators. The matrix $\Delta$ will account for the interactions, and will hence be of higher-order in $\eta$. Our approximant will be of the form $\widetilde{H}+\Delta$. 

Our next step to determine these matrices will be to compute the eigenvalues of the full Hamiltonian $H^0$.

\subsection{Two-resonators case}\label{sec:tightbind1}
Firstly, we compute the eigenvalues of the full static Hamiltonian \eqref{static_H}. We compute its characteristic polynomial. By noticing that $H^0 - \lambda \mathbb{Id}_4$ can be viewed as a $2 \times 2$ block matrix whose blocks commute with each other, we may compute the block version of the determinant and turn the eigenvalue problem for the Hamiltonian $H^0$ into the one for $M^0$:
\begin{equation}\label{reducing characteristic of H full}
    0= \det (- \lambda \mathbb{Id}_{4} + H^0) = \begin{vmatrix} -\lambda \mathbb{Id}_2 & \sqrt{M^0} \\
    \sqrt{M^0} & -\lambda \mathbb{Id}_2
    \end{vmatrix} = \det \left ( \lambda^2 \mathbb{Id}_2 - M^0 \right ).
\end{equation}

To compute the eigenvalues of $M^0$, we will employ (as in the rest of this section) the diluteness assumption \eqref{capacitance matrix in dilute regime}. Matrices of the form 
$$
    \begin{pmatrix}
    \alpha & \beta \\
    \beta & \alpha 
    \end{pmatrix},
$$
where $\alpha$, $\beta \in \mathbb{C}$, have $\alpha \pm \beta$ as two eigenvalues. From \eqref{capacitance matrix in dilute regime}, $M^0$ has this form and therefore we can use this to compute its eigenvalues and  find that

\begin{equation}\label{Eigenvalues full M^0 2 resonators}
\lambda^2 = \delta \frac{\kappa_r \capacity_B }{\rho_r |B|}\pm \delta \eta \frac{\kappa_r \left ( \capacity_B \right )^2}{4 \pi \rho_r |z_1-z_2| |B|} + O(\eta^2),
\end{equation}
where $z_1$ and $z_2$ are the positions of the centres of the two resonators. 
To keep track of the four eigenvalues, we define 
\begin{align}\label{eigenvalues full Hamiltonian 2 res}
\begin{cases}
    &\ds \lambda_{11}:= \sqrt{\delta \frac{\kappa_r \capacity_B}{\rho_r |B|} }\left(1+ \eta\frac{\capacity_B}{8 \pi |z_1-z_2|} +O(\eta^2) \right),  \\
    \nm 
    &\ds \lambda_{12}:= \sqrt{\delta \frac{\kappa_r \capacity_B}{\rho_r |B|} }\left(1- \eta \frac{\capacity_B}{8 \pi |z_1-z_2|} +O(\eta^2) \right), \\
    \nm
    &\ds \lambda_{21}:= - \sqrt{\delta \frac{\kappa_r \capacity_B}{\rho_r |B|} }\left(1+ \eta \frac{\capacity_B}{8 \pi |z_1-z_2|} +O(\eta^2) \right), \\
    \nm
    &\ds \lambda_{22}:= -\sqrt{\delta \frac{\kappa_r \capacity_B}{\rho_r |B|} }\left(1- \eta \frac{\capacity_B}{8 \pi |z_1-z_2|} +O(\eta^2) \right) .
  \end{cases}
\end{align}

Following \cite{ammari2017plasmonicscalar,pnas}, we intend to find a tight-binding approximation for the full static Hamiltonian. The choice that proved to fit better with the criteria listed at the beginning of this section is the following:
\begin{align}\label{Tilde H}
    \widetilde{H} := \sqrt{\delta \frac{\kappa_r \capacity_B}{\rho_r |B|} }\text{ diag}\begin{pmatrix}
    1, &1, &-1, &-1
    \end{pmatrix}.
\end{align}
It is clear from  formula \eqref{eigenvalues full Hamiltonian 2 res} that the eigenvalues of $\widetilde{H}^0$ coincide with the ones of $H^0$ at zero-th order in $\eta$. Hence, the first of the criteria listed at the beginning of this section is satisfied. 

We would now like to find a $4\times 4$ matrix $\Delta$ such that $\widetilde{H}^0 + \Delta$ has, up to first-order in $\eta$, the same eigenvalues as $H^0$. In line with the three requirements listed at the beginning of this section, our choice of $\Delta$ is 
\begin{equation*}
    \Delta = \begin{pmatrix}
    \Delta_1 & \mathbf{0}\\ \mathbf{0} & \Delta_2
    \end{pmatrix},
\end{equation*}
where $\Delta_1$ and $\Delta_2$ are two $2 \times 2$ matrices of the form
\begin{align}
    \Delta_1 = \begin{pmatrix}
    0 & a\\
    a & 0
    \end{pmatrix}, & \qquad
    \Delta_2 = \begin{pmatrix}
    0 & b\\
    b & 0
    \end{pmatrix}
\end{align}
with $a,$ $b \in \mathbb{R}$ being two parameters to be determined. This way, $\Delta$ represents the interactions between the resonators and does not couple the negative modes with the positive ones. In this particular case of two resonators, the matrix $\widetilde{H}+\Delta$ is tridiagonal. A nearest neighbour approximation holds trivially in this case. The parameters $a$ and $b$ can be viewed as the hopping parameters of the system.

To determine $a$ and $b$, let us look at the eigenvalues of $\widetilde{H}+ \Delta$. We have
\begin{align*}
    0 & = \det \left( \widetilde{H} + \Delta - \widetilde{\lambda} \mathbb{Id}_4\right)
      \\& = \begin{vmatrix} \sqrt{\delta \frac{\kappa_r \capacity_B}{\rho_r |B|} } - \widetilde{\lambda} & a 
    \\ a & \sqrt{\delta \frac{\kappa_r \capacity_B}{\rho_r |B|} } - \widetilde{\lambda} \end{vmatrix}\begin{vmatrix} -\sqrt{\delta \frac{\kappa_r \capacity_B}{\rho_r |B|} } - \widetilde{\lambda} & b 
    \\ b & -\sqrt{\delta \frac{\kappa_r \capacity_B}{\rho_r |B|} } - \widetilde{\lambda} \end{vmatrix}  \\
    & = \left[\left( \sqrt{\delta \frac{\kappa_r \capacity_B}{\rho_r |B|} } - \widetilde{\lambda}\right)^2 - a^2\right] \left[\left(- \sqrt{\delta \frac{\kappa_r \capacity_B}{\rho_r |B|} } - \widetilde{\lambda} \right)^2 - b^2 \right],
\end{align*}
from which we obtain
\begin{equation}\label{approximate first order eigenvalues}
\begin{aligned}
   & \widetilde{\lambda}_{11}:= \sqrt{\delta \frac{\kappa_r \capacity_B}{\rho_r |B|} } + a ,
    & \widetilde{\lambda}_{12}:= \sqrt{\delta \frac{\kappa_r \capacity_B}{\rho_r |B|} } - a , \\
    & \widetilde{\lambda}_{21}:= - \sqrt{\delta \frac{\kappa_r \capacity_B}{\rho_r |B|} } - b ,
    & \widetilde{\lambda}_{22}:=  - \sqrt{\delta \frac{\kappa_r \capacity_B}{\rho_r |B|} } + b  .
\end{aligned}
\end{equation}
By comparing the approximate eigenvalues \eqref{approximate first order eigenvalues} to the eigenvalues of the full matrix $H^0$ given in \eqref{eigenvalues full Hamiltonian 2 res} we obtain
\begin{align*}
	a = \eta \frac{\capacity_B}{8 \pi |z_1-z_2|}, \qquad b = \eta \frac{\capacity_B}{8 \pi |z_1-z_2|}.
\end{align*}

Plugging the values of $a$ and $b$ in the formulas \eqref{approximate first order eigenvalues}, we obtain
\begin{equation}\label{final approximate first order eigenvalues}
\begin{aligned}
    & \widetilde{\lambda}_{11}:= \sqrt{\delta \frac{\kappa_r \capacity_B}{\rho_r |B|} } + \eta \frac{\capacity_B}{8 \pi |z_1-z_2|},
    & \widetilde{\lambda}_{12}:=  \sqrt{\delta \frac{\kappa_r \capacity_B}{\rho_r |B|} } - \eta \frac{\capacity_B}{8 \pi |z_1-z_2|},  \\
    & \widetilde{\lambda}_{21}:= - \sqrt{\delta \frac{\kappa_r \capacity_B}{\rho_r |B|} } - \eta \frac{\capacity_B}{8 \pi |z_1-z_2|},
    & \widetilde{\lambda}_{22}:=  - \sqrt{\delta \frac{\kappa_r \capacity_B}{\rho_r |B|} } + \eta \frac{\capacity_B}{8 \pi |z_1-z_2|}  .
\end{aligned}
\end{equation}

Lastly, the approximant matrix reads
\begin{equation*}
    \widetilde{H}+\Delta = \begin{pmatrix}
    \sqrt{\delta \frac{\kappa_r \capacity_B}{\rho_r |B|} }&\eta \frac{\capacity_B}{8 \pi |z_1-z_2|} & 0 &0 \\
    \eta \frac{\capacity_B}{8 \pi |z_1-z_2|} & \sqrt{\delta \frac{\kappa_r \capacity_B}{\rho_r |B|} }& 0&0 \\
    0&0&-\sqrt{\delta \frac{\kappa_r \capacity_B}{\rho_r |B|} }& \eta \frac{\capacity_B}{8 \pi |z_1-z_2|}\\
    0&0& \eta \frac{\capacity_B}{8 \pi |z_1-z_2|} &-\sqrt{\delta \frac{\kappa_r \capacity_B}{\rho_r |B|} }
    \end{pmatrix}.
\end{equation*}
    
Let us remark that the approximant matrix, up to second-order in $\eta$, is similar to the starting Hamiltonian.

\subsection{Case of an arbitrary number of resonators}
In this section we will outline the case of an arbitrary number of resonators. We will proceed analogously to the two-resonator case. To do so, we will first compute the eigenvalues of the static Hamiltonian matrix $H^0$ up to first-order in the dilute regime parameter $\eta$. We first notice that
\begin{equation}
    0= \det (- \lambda \mathbb{Id}_{2N} + H^0) = \begin{vmatrix} -\lambda \mathbb{Id}_N & \sqrt{M^0} \\
    \sqrt{M^0} & -\lambda \mathbb{Id}_N
    \end{vmatrix} = \det \left ( \lambda^2 \mathbb{Id}_N - M^0 \right ),
\end{equation}
since the submatrices shown in the computation all commute with each other. Therefore, the problem of computing the eigenvalues of the Hamiltonian is equivalent to the one of computing the eigenvalues of the generalised capacitance matrix. Let us therefore compute the eigenvalues of the generalised capacitance matrix.

It follows from \eqref{capacitance matrix in dilute regime} that $M^0$ is of the form
$$
    M^0= \delta \frac{\kappa_r\capacity_B}{\rho_r |B|} \left(\mathbb{Id}_N + \eta A \right) + O(\eta^2),
$$
where $A$ is the $N \times N$ matrix whose $(j,k)$-th entry is 
$$
    A_{jk} = \begin{cases}
        \frac{\capacity_B}{4\pi |z_j - z_k|} & \text{if } j\not =k, \\
        0 &  \text{if } j =k.
    \end{cases}
$$
As $M^0$ is diagonalisable, it follows that the eigenvalues of $M^0$ are
$$
    \lambda^2_j = \delta \frac{\kappa_r\capacity_B}{\rho_r |B|} \left(1 + \eta a_j \right) + O(\eta^2),
$$
where $1\leq j \leq N$ and $a_j$ is the $j$-th eigenvalue of $A$. Consequently, the eigenvalues of $H^0$ are
\begin{equation}\label{eq:lambdaj}
\lambda_j = \pm\sqrt{\delta \frac{\kappa_r \capacity_B}{\rho_r |B|} }\left(1 + \eta \frac{a_j}{2} \right) + O(\eta^2).
\end{equation}

To choose an approximant, as in the case for $N=2$, we define
$$
    \widetilde{H} := \sqrt{\delta \frac{\kappa_r \capacity_B}{\rho_r |B|} } \text{ diag} (\underbrace{1, \ 1,\ ..., \ 1}_N, \underbrace{-1, \ -1,\ ...,\ -1}_N).
$$
Analogously as in \Cref{sec:tightbind1}, the higher-order coupling matrix $\Delta$ will in general be $2 \times 2$ block diagonal, so that no coupling occurs between the negative modes and the positive ones. We therefore seek $\Delta$ of the form
$$\Delta = \begin{pmatrix}
	\Delta_1 & \mathbf{0}\\ \mathbf{0} & \Delta_2
\end{pmatrix},$$
where $\Delta_1, \Delta_2$ are symmetric, off-diagonal matrices. Comparing the eigenvalues of  $\widetilde{H} + \Delta$ to $\lambda_j$, given in \eqref{eq:lambdaj}, we obtain that (up to similarity transformations)
$$\Delta_1 = \Delta_2 = \frac{\eta}{2}\sqrt{\delta \frac{\kappa_r \capacity_B}{\rho_r |B|} }A.$$
It is evident that the two diagonal blocks $\Delta_1, \Delta_2$ are dense and the tight-binding approximant is not tridiagonal. In other words, long-range interactions cannot be ignored, and nearest-neighbour approximations are not appropriate. 

\begin{rmk}
In the case where the material parameters are modulated inside the resonators, we would proceed in a slightly different way. We would diagonalise the modulated Hamiltonian $H^0+\varepsilon H^1$ first. The leading-order in both the dilute regime parameter $\eta$ and the modulation parameter $\varepsilon$ (which will depend upon each other) will provide the sign of the eigenvalues, which will be ordered in $\widetilde{H}$ as in the static case (the positive first, then the negative). We refer to \cite{TheaThesis} for a glimpse at the computations needed to obtain these eigenvalues. For the coupling matrix $\Delta$, which now will in general depend on time, we would proceed analogously as before: we impose it to be $2 \times 2$ block diagonal, and we solve the systems obtained by imposing the equality between the eigenvalues. By splitting the equations for the static part and the modulated part, we would obtain two systems. The static part which has already been computed should solve the first, while the modulated part could be obtained by solving the second system. We expect it to be as complicated as the system \eqref{General Component-wise}, since in the end the matrix $H$ and the approximant will only differ by a similarity transform.
\end{rmk}

\section{Concluding remarks}
In the first part of this paper, we recast the capacitance matrix formulation of the Hill equations \eqref{eq:hill} arising from the subwavelength problem \eqref{eq:wave} in a Hamiltonian formulation \eqref{MatrixDefH}. To do so, we had to determine a proper transformation $T$ that would make $H$ Hermitian. The system arising from imposing the Hermitianity conditions is rather complicated to solve. We used a perturbative argument to simplify the problem by assuming the modulation to be asymptotically small. The static system was completely solved by \eqref{static_H}, where an extra physical interpretation could be obtained by noticing its relationship with the generalised capacitance matrix. For the modulated part, the system still proved hard to solve. We obtained explicit closed formulas for the one-resonator case. For the case of arbitrary number of resonators, we provided conditions that made the under-determined system determined, which proved existence of a transformation yielding a Hermitian Hamiltonian \eqref{MatrixDefH}. For a fixed number of resonators, such equations can be integrated to obtain an explicit analytic solution, although its complexity makes it rather impractical in general.

In the second part of this paper, we have highlighted similarities and differences between the capacitance matrix formulation for the subwavelength resonance problem \cite{ammari2021functional} and the tight-binding regime in condensed matter theory. The main difference is that linear combinations of modes of the subwavelength resonance problem (with no modulation) cannot, in general, be a good approximation of modes relative to multiple resonators, as they need to be almost constant inside each of the resonators (see \cite{ammari2021functional}). Therefore, we needed to assume the resonators to be asymptotically small in size (the so called dilute regime given by \eqref{capacitance matrix in dilute regime}). Starting from the Hamiltonian formulation previously discovered, we created a model based on the tight-binding regime. We achieved this using a perturbative argument in terms of the diluteness parameter. The entries of the diagonal of our approximant represent the self-interaction between a single resonator, which also happens to be the lowest-order term of the eigenvalues of our Hamiltonian. Subsequently, we added a higher-order $2 \times 2$ block diagonal matrix which represents the interactions and couples the modes. The off-diagonal blocks are chosen to be zero in order to avoid  non-physical coupling between modes at positive and at negative frequencies. We emphasised that the non-zero blocks of this matrix are dense due to the fact that the capacitance matrix is dense. Therefore, long-range interaction cannot be ignored, and the nearest neighbour approximation is not applicable.
\bibliographystyle{plain}
\bibliography{references}
\end{document}

%% file: resonator.tikzstyles
\tikzstyle{resonator}=[-, draw=none, fill={rgb,255: red,191; green,191; blue,191}]
\tikzstyle{mode1}=[-, draw=blue]

%% file: figures/Single-resonator_modes.tikz
\begin{tikzpicture}
	\begin{pgfonlayer}{nodelayer}
		\node [style=none] (0) at (-3, 0) {};
		\node [style=none] (1) at (-1, 0) {};
		\node [style=none] (2) at (-2, 0.5) {};
		\node [style=none] (3) at (-2, -0.5) {};
		\node [style=none] (4) at (1, 0) {};
		\node [style=none] (5) at (3, 0) {};
		\node [style=none] (6) at (2, 0.5) {};
		\node [style=none] (7) at (2, -0.5) {};
		\node [style=none] (9) at (-2.5, 4) {};
		\node [style=none] (10) at (-7, 0.25) {};
		\node [style=none] (12) at (-1.5, 4) {};
		\node [style=none] (15) at (3, 0.25) {};
		\node [style=none] (16) at (-9.75, 0.1) {};
		\node [style=none] (17) at (9.75, 0.075) {};
		\node [style=none] (18) at (1.5, 4) {};
		\node [style=none] (19) at (-3, 0.25) {};
		\node [style=none] (20) at (2.5, 4) {};
		\node [style=none] (21) at (7, 0.25) {};
		\node [style=none] (22) at (-9.5, 0.1) {};
		\node [style=none] (23) at (9.75, 0.075) {};
	\end{pgfonlayer}
	\begin{pgfonlayer}{edgelayer}
		\draw [style=resonator] (0.center)
			 to [in=-180, out=-90, looseness=0.75] (3.center)
			 to [in=-90, out=0, looseness=0.75] (1.center)
			 to [in=0, out=90, looseness=0.75] (2.center)
			 to [in=90, out=180, looseness=0.75] cycle;
		\draw [style=resonator] (4.center)
			 to [in=-180, out=-90, looseness=0.75] (7.center)
			 to [in=-90, out=0, looseness=0.75] (5.center)
			 to [in=0, out=90, looseness=0.75] (6.center)
			 to [in=90, out=180, looseness=0.75] cycle;
		\draw (16.center)
			 to (10.center)
			 to [in=-180, out=0, looseness=0.75] (9.center)
			 to (12.center)
			 to [in=-180, out=0, looseness=0.75] (15.center)
			 to (17.center);
		\draw [style=mode1] (22.center)
			 to (19.center)
			 to [in=-180, out=0, looseness=0.75] (18.center)
			 to (20.center)
			 to [in=-180, out=0, looseness=0.75] (21.center)
			 to (23.center);
	\end{pgfonlayer}
\end{tikzpicture}

%% file: figures/Two-resonator_mode.tikz
\begin{tikzpicture}
	\begin{pgfonlayer}{nodelayer}
		\node [style=none] (0) at (-3, 0) {};
		\node [style=none] (1) at (-1, 0) {};
		\node [style=none] (2) at (-2, 0.5) {};
		\node [style=none] (3) at (-2, -0.5) {};
		\node [style=none] (4) at (1, 0) {};
		\node [style=none] (5) at (3, 0) {};
		\node [style=none] (6) at (2, 0.5) {};
		\node [style=none] (7) at (2, -0.5) {};
		\node [style=none] (9) at (-2.5, 4) {};
		\node [style=none] (10) at (-7, 0.25) {};
		\node [style=none] (11) at (7, 0.25) {};
		\node [style=none] (12) at (-1.25, 4) {};
		\node [style=none] (13) at (1.25, 4) {};
		\node [style=none] (14) at (2.5, 4) {};
		\node [style=none] (15) at (0, 3) {};
	\end{pgfonlayer}
	\begin{pgfonlayer}{edgelayer}
		\draw [style=resonator] (0.center)
			 to [in=-180, out=-90, looseness=0.75] (3.center)
			 to [in=-90, out=0, looseness=0.75] (1.center)
			 to [in=0, out=90, looseness=0.75] (2.center)
			 to [in=90, out=180, looseness=0.75] cycle;
		\draw [style=resonator] (4.center)
			 to [in=-180, out=-90, looseness=0.75] (7.center)
			 to [in=-90, out=0, looseness=0.75] (5.center)
			 to [in=0, out=90, looseness=0.75] (6.center)
			 to [in=90, out=180, looseness=0.75] cycle;
		\draw (10.center)
			 to [in=-165, out=0, looseness=0.75] (9.center)
			 to (12.center)
			 to [in=-180, out=0, looseness=0.75] (15.center)
			 to [in=-180, out=0, looseness=0.75] (13.center)
			 to (14.center)
			 to [in=-180, out=-15, looseness=0.75] (11.center);
	\end{pgfonlayer}
\end{tikzpicture}

%% file: figures/Linear_combination.tikz
\begin{tikzpicture}
	\begin{pgfonlayer}{nodelayer}
		\node [style=none] (0) at (-3, 0) {};
		\node [style=none] (1) at (-1, 0) {};
		\node [style=none] (2) at (-2, 0.5) {};
		\node [style=none] (3) at (-2, -0.5) {};
		\node [style=none] (4) at (1, 0) {};
		\node [style=none] (5) at (3, 0) {};
		\node [style=none] (6) at (2, 0.5) {};
		\node [style=none] (7) at (2, -0.5) {};
		\node [style=none] (18) at (0, 3.5) {};
		\node [style=none] (23) at (7, 0.075) {};
		\node [style=none] (24) at (1, 4) {};
		\node [style=none] (25) at (-1, 4) {};
		\node [style=none] (26) at (-7, 0.075) {};
	\end{pgfonlayer}
	\begin{pgfonlayer}{edgelayer}
		\draw [style=resonator] (0.center)
			 to [in=-180, out=-90, looseness=0.75] (3.center)
			 to [in=-90, out=0, looseness=0.75] (1.center)
			 to [in=0, out=90, looseness=0.75] (2.center)
			 to [in=90, out=180, looseness=0.75] cycle;
		\draw [style=resonator] (4.center)
			 to [in=-180, out=-90, looseness=0.75] (7.center)
			 to [in=-90, out=0, looseness=0.75] (5.center)
			 to [in=0, out=90, looseness=0.75] (6.center)
			 to [in=90, out=180, looseness=0.75] cycle;
		\draw [in=-165, out=0, looseness=0.75] (26.center) to (25.center);
		\draw [in=180, out=0, looseness=0.75] (25.center) to (18.center);
		\draw [in=-180, out=0] (18.center) to (24.center);
		\draw [in=180, out=-15, looseness=0.75] (24.center) to (23.center);
	\end{pgfonlayer}
\end{tikzpicture}

%% file: main_final.bbl
\begin{thebibliography}{10}

\bibitem{ammari2019cochlea}
Habib Ammari and Bryn Davies.
\newblock A fully coupled subwavelength resonance approach to filtering
  auditory signals.
\newblock {\em Proc. R. Soc. A}, 475(2228):20190049, 2019.

\bibitem{ammari2020biomimetic}
Habib Ammari and Bryn Davies.
\newblock A biomimetic basis for auditory processing and the perception of
  natural sounds.
\newblock {\em arXiv preprint arXiv:2005.12794}, 2020.

\bibitem{ammari2020robust}
Habib Ammari, Bryn Davies, and Erik~Orvehed Hiltunen.
\newblock Robust edge modes in dislocated systems of subwavelength resonators.
\newblock {\em arXiv preprint arXiv:2001.10455}, 2020.

\bibitem{ammari2021functional}
Habib Ammari, Bryn Davies, and Erik~Orvehed Hiltunen.
\newblock Functional analytic methods for discrete approximations of
  subwavelength resonator systems, 2021.

\bibitem{ammari2020high}
Habib Ammari, Bryn Davies, Erik~Orvehed Hiltunen, Hyundae Lee, and Sanghyeon
  Yu.
\newblock High-order exceptional points and enhanced sensing in subwavelength
  resonator arrays.
\newblock {\em Stud. Appl. Math.}, 146(2):440--462, 2021.

\bibitem{ammari2020topological}
Habib Ammari, Bryn Davies, Erik~Orvehed Hiltunen, and Sanghyeon Yu.
\newblock Topologically protected edge modes in one-dimensional chains of
  subwavelength resonators.
\newblock {\em J. Math. Pures Appl.}, 144:17--49, 2020.

\bibitem{ammari2020exceptional}
Habib Ammari, Bryn Davies, Hyundae Lee, Erik~Orvehed Hiltunen, and Sanghyeon
  Yu.
\newblock Exceptional points in parity--time-symmetric subwavelength
  metamaterials.
\newblock {\em arXiv preprint arXiv:2003.07796}, 2020.

\bibitem{ammari2020honeycomb}
Habib Ammari, Brian Fitzpatrick, Erik~Orvehed Hiltunen, Hyundae Lee, and
  Sanghyeon Yu.
\newblock Honeycomb-lattice {M}innaert bubbles.
\newblock {\em SIAM J. Math. Anal.}, 52(6):5441--5466, 2020.

\bibitem{ammari2018mathematical}
Habib Ammari, Brian Fitzpatrick, Hyeonbae Kang, Matias Ruiz, Sanghyeon Yu, and
  Hai Zhang.
\newblock {\em Mathematical and Computational Methods in Photonics and
  Phononics}, volume 235 of {\em Mathematical Surveys and Monographs}.
\newblock American Mathematical Society, Providence, 2018.

\bibitem{ammari2019double}
Habib Ammari, Brian Fitzpatrick, Hyundae Lee, Sanghyeon Yu, and Hai Zhang.
\newblock Double-negative acoustic metamaterials.
\newblock {\em Quart. Appl. Math.}, 77(4):767--791, 2019.

\bibitem{ammari2020time}
Habib Ammari and Erik~Orvehed Hiltunen.
\newblock Time-dependent high-contrast subwavelength resonators.
\newblock {\em J. Comp. Phys., to appear (arXiv preprint arXiv:2012.10274)},
  2020.

\bibitem{TheaThesis}
Habib Ammari, Erik~Orvehed Hiltunen, and Thea Kosche.
\newblock Asymptotic floquet theory for first order odes with finite fourier
  series perturbation and applications in time-modulated metamaterials.
\newblock {\em preprint}, 2021.

\bibitem{ammari2018high}
Habib Ammari, Erik~Orvehed Hiltunen, and Sanghyeon Yu.
\newblock A high-frequency homogenization approach near the {Dirac} points in
  bubbly honeycomb crystals.
\newblock {\em Arch. Rational Mech. Anal.}, 238:1559--1583, 2020.

\bibitem{ammari2019bloch}
Habib Ammari, Hyundae Lee, and Hai Zhang.
\newblock Bloch waves in bubbly crystal near the first band gap: a
  high-frequency homogenization approach.
\newblock {\em SIAM J. Math. Anal.}, 51(1):45--59, 2019.

\bibitem{ammari2017plasmonicscalar}
Habib Ammari, Pierre Millien, Matias Ruiz, and Hai Zhang.
\newblock Mathematical analysis of plasmonic nanoparticles: the scalar case.
\newblock {\em Arch. Rational Mech. Anal.}, 224(2):597--658, 2017.

\bibitem{ammari2017effective}
Habib Ammari and Hai Zhang.
\newblock Effective medium theory for acoustic waves in bubbly fluids near
  {Minnaert} resonant frequency.
\newblock {\em SIAM J. Math. Anal.}, 49(4):3252--3276, 2017.

\bibitem{cassier21}
M.~Cassier and M.~I. Weinstein.
\newblock High contrast elliptic operators in honeycomb structures.
\newblock {\em arXiv:2103.16682}, 2021.

\bibitem{https://doi.org/10.1002/cpa.21735}
Charles~L. Fefferman, James~P. Lee-Thorp, and Michael~I. Weinstein.
\newblock Honeycomb schrödinger operators in the strong binding regime.
\newblock {\em Communications on Pure and Applied Mathematics},
  71(6):1178--1270, 2018.

\bibitem{review2}
Hao Ge, Min Yang, Chu Ma, Ming-Hui Lu, Yan-Feng Chen, Nicholas Fang, and Ping
  Sheng.
\newblock Breaking the barriers: advances in acoustic functional materials.
\newblock {\em National Science Review}, 5:159--182, 2018.

\bibitem{gohberg1971operator}
I.C. Gohberg and E.I. Sigal.
\newblock An operator generalization of the logarithmic residue theorem and the
  theorem of {R}ouch\'{e}.
\newblock {\em Sb. Math.}, 13(4):603--625, 1971.

\bibitem{lemoult2016soda}
Fabrice Lemoult, Nadège Kaina, Mathias Fink, and Geoffroy Lerosey.
\newblock Soda cans metamaterial: A subwavelength-scaled phononic crystal.
\newblock {\em Crystals}, 6(7), 2016.

\bibitem{phononic2}
Z.~Liu, C.T. Chan, and P.~Sheng.
\newblock Analytic model of phononic crystals with local resonances.
\newblock {\em Phys. review B}, 71(1):014103, 2005.

\bibitem{phononic1}
Z.~Liu, X.~Zhang, Y.~Mao, Y.Y. Zhu, Z.~Yang, C.T. Chan, and P.~Sheng.
\newblock Locally resonant sonic materials.
\newblock {\em Science}, 289(5485):1734--1736, 2000.

\bibitem{review}
Guancong Ma and Ping Sheng.
\newblock Acoustic metamaterials: From local resonances to broad horizons.
\newblock {\em Science Advances}, 2(2):e1501595, 2016.

\bibitem{minnaert1933musical}
M.~Minnaert.
\newblock {O}n musical air-bubbles and the sounds of running water.
\newblock {\em Philos. Mag.}, 16(104):235--248, 1933.

\bibitem{osti_377103}
A~L Mironov and V~L Olenik.
\newblock Limits of applicability of the tight binding approximation.
\newblock {\em Theoretical and Mathematical Physics}, 99(1), 10 1994.

\bibitem{teschl2012ordinary}
G.~Teschl.
\newblock {\em Ordinary Differential Equations and Dynamical Systems}.
\newblock Graduate studies in mathematics. American Mathematical Society, 2012.

\bibitem{pnas}
Sanghyeon Yu and Habib Ammari.
\newblock Hybridization of singular plasmons via transformation optics.
\newblock {\em Proc. Natl. Acad. Sci. USA}, 116:13785--13790, 2019.

\bibitem{yves2017crytalline}
Simon Yves, Romain Fleury, Thomas Berthelot, Mathias Fink, Fabrice Lemoult, and
  Geoffroy Lerosey.
\newblock Crystalline metamaterials for topological properties at subwavelength
  scales.
\newblock {\em Nat. Commun.}, 8(1):16023, Jul 2017.

\end{thebibliography}
